\documentclass[11pt, english]{article}

\usepackage{babel}
\usepackage[latin1]{inputenc}
\usepackage{amsmath}
\usepackage{amssymb}
\usepackage{amsthm}
\usepackage{graphicx}

\theoremstyle{plain}
\newtheorem{Thm}{Theorem}%[section]
\newtheorem{Prop}[Thm]{Proposition}

\newtheorem{Lemma}[Thm]{Lemma}

\def\N{\mathbb{N}}
\def\n{\mathbb{N}_0}
\def\G{\text{-GDWN}}
\usepackage{pgf}
\usepackage{tikz}
\usetikzlibrary{arrows}
\usetikzlibrary{decorations.pathmorphing}
\usetikzlibrary{backgrounds}
\usetikzlibrary{fit}

\title{$(1,2)$-GDWN splits}
\author{Urban Larsson\\ Mathematical Sciences,\\ Chalmers University of Technology and University of Gothenburg,\\ G\"oteborg, Sweden\\
urban.larsson@chalmers.se}
\date{\today}
\begin{document}
\maketitle
\begin{abstract}
 We study impartial take away games on 2 unordered piles of finite nonnegative numbers of tokens $(x,y)$. Two players alternate in removing at least one and at most all tokens from the respective piles, according to certain rules, and the game terminates when a player in turn is unable to move. We follow the normal play convention, which means that a player who cannot move loses. In the game of Wythoff Nim, a player is allowed to remove either any number of tokens from precisely one of the piles or the same number of tokens from both. Let $\phi = \frac{1+\sqrt{5}}{2}$ and for all nonnegative integers $n$, $A_n=\lfloor\phi n \rfloor$ and $B_n=A_n+n$. The P-positions of Wythoff Nim are all pairs of piles with $A_n$ and $B_n$ tokens respectively. We study a generalization of this game called $(1,2)\G$ where, in addition to the rules of Wythoff Nim, a player has the choice to remove a positive number of tokens from one of the piles and twice that number from the other pile. We show that there is an infinite sector $\alpha \le y/x \le \alpha +\epsilon$, for given real numbers $\alpha>1$ and $\epsilon > 0$, for which each $(x,y)$ is an N-position, but that there are infinitely many P-positions for both $1\le y/x <\alpha $ and  $\alpha +\epsilon < y/x $. This proves a conjecture from a recent paper. Namely, the adjoined set of moves in $(1,2)\G$ \emph{splits} the beam of slope $\phi $  P-positions of Wythoff Nim (in the same sense that the adjoined moves in Wythoff Nim split the single slope 1 beam of P-positions of 2-pile Nim).
We also provide a lower bound on the lower asymtotic density of lower pile heights of P-positions for extensions of Wythoff Nim.  Suppose that $(a_i)$ and $(b_i)$, $i>0$, is a pair of so-called complementary sequences on the natural numbers which satisfy $(a_i)$ is increasing and for all $i$, $a_i<b_i$, for all $i\ne j$, $b_i-a_i\ne b_j-a_j$. Then $\liminf_{n\rightarrow \infty}\frac{\#\{i\mid a_i < n\}}{n} \ge \phi^{-1}$.

\end{abstract}

%\maketitle
\section{Introduction}
We study generalizations of the 2-player impartial take-away games of \emph{2-pile Nim} \cite{Bou02} and \emph{Wythoff Nim}, \cite{Wyt07, HeLa06, Lar09, Lar}. A background on impartial (take-away) games can be found in for example \cite{BCG82, Con76}. We use some standard terminology for such games without ties. A position is a previous-player win, a \emph{P-position}, if none of its options are P-positions; otherwise it is a next-player win, an \emph{N-position}. We follow the conventions of \emph{normal play}, that is, a player who is not able to move loses and the other player wins. Thus, given an impartial game, we get a recursive characterization of the set of all P-positions beginning with the terminal position(s). 

\begin{figure}[ht!]
\begin{center}
\vspace{0.5 cm}
{\includegraphics[width=0.75\textwidth]{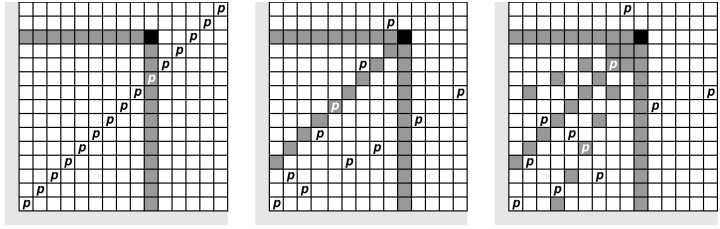}}
\end{center}\caption{The figures illustrate typical 
moves (in dark gray) and initial P~-~positions of Nim, Wythoff Nim, and $(1, 2)$-GDWN respectively. The black square is a given game position. The white P's represent the winning options from this position for the respective games. Hence the given position is in N for each game, but for different reasons. In Nim, a certain Nim-type move suffices, whereas in Wythoff Nim a certain diagonal type move is required. Similarly, in GDWN either a slope 2 type move or a slope 1 type move suffices.}\label{F0}
\end{figure}

\begin{figure}[ht!]
\begin{center}
{\includegraphics[width=0.19\textwidth]{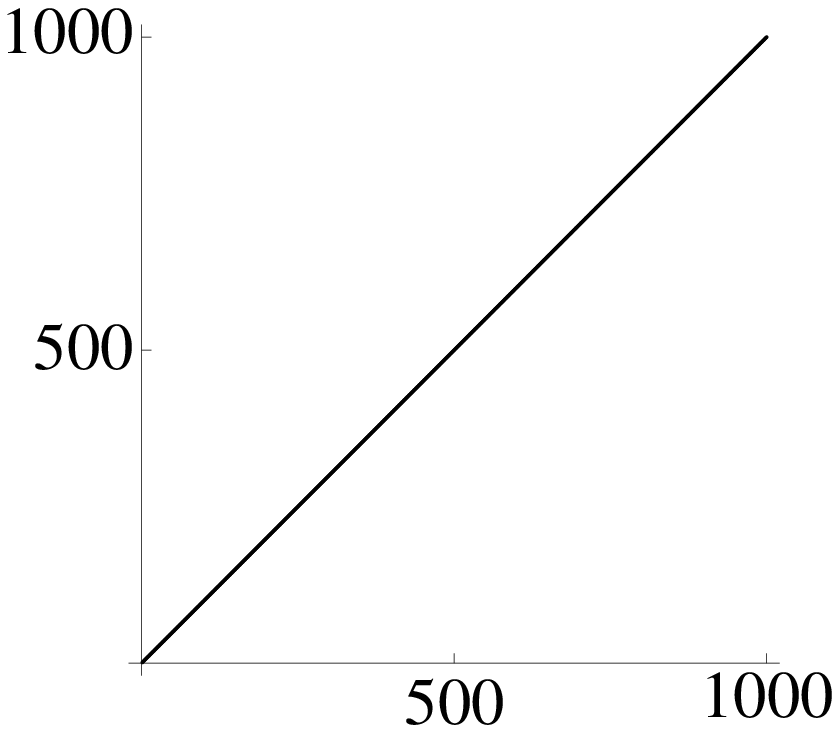}\hspace{0.4 cm}
{\includegraphics[width=0.19\textwidth]{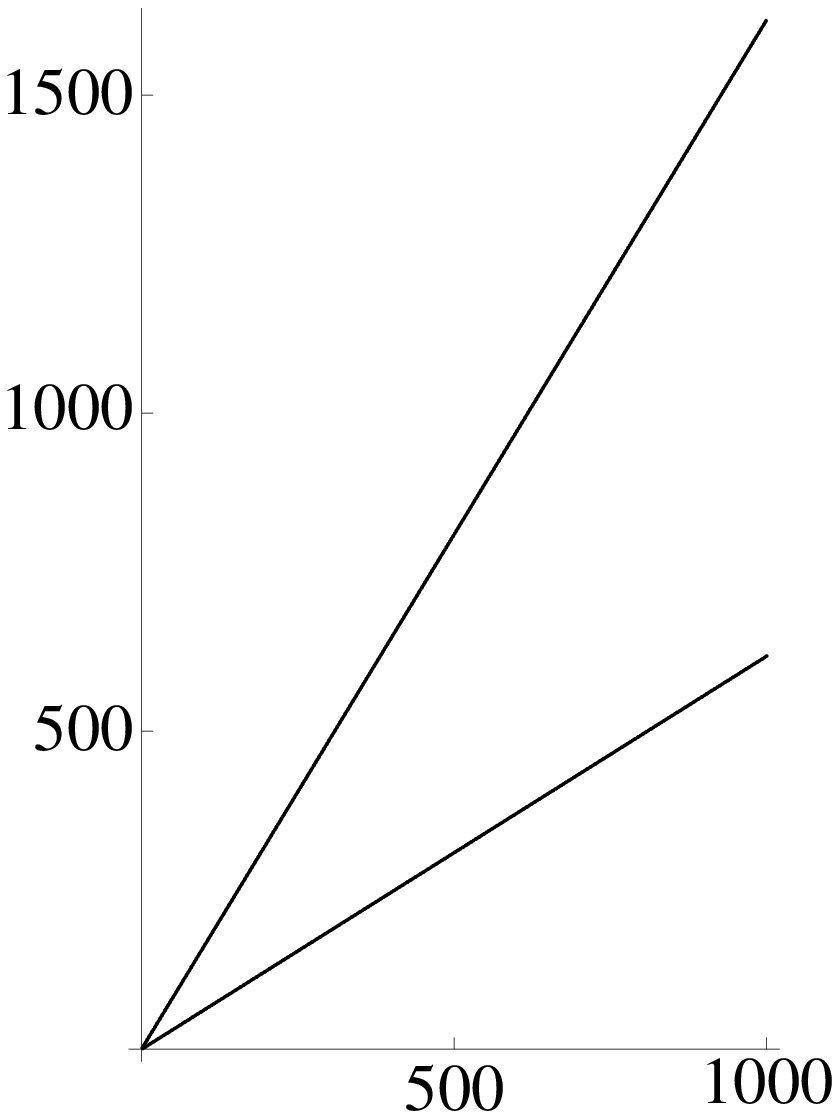}}
\includegraphics[width=0.32\textwidth]{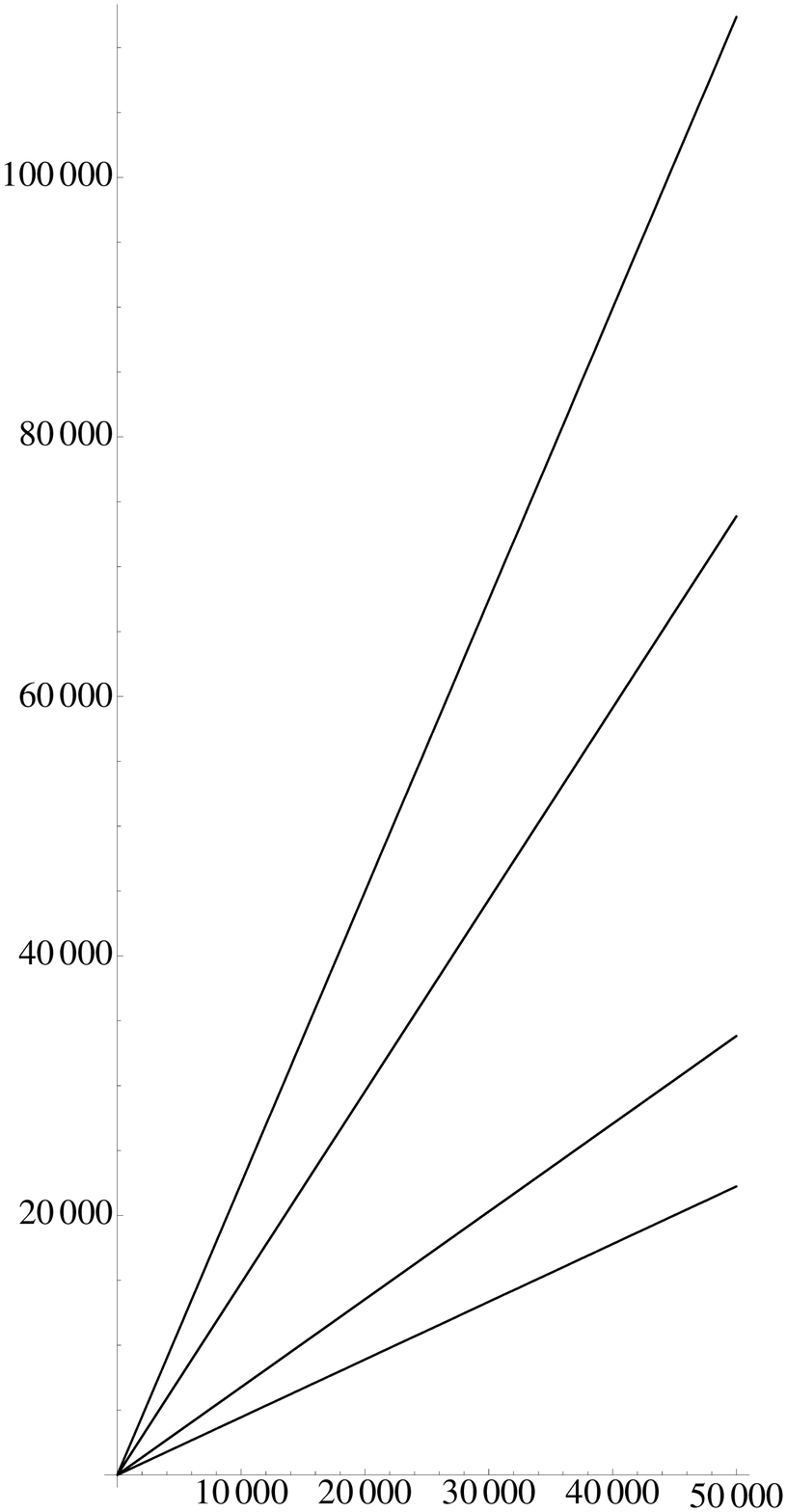}}
\end{center}\caption{These figures give the initial P-positions of the games Nim, Wythoff Nim and $(1,2)\G$. The left most figure illustrates 2-pile Nim's single P-beam of slope 1. Then, in the middle we illustrate Wythoff Nim's pair of P-beams with slopes $\phi^{-1}$ and $\phi$ respectively and, at last, we present the initial P-positions of $(1,2)\G$, for all $x$-coordinates $\le 50000$, with its split of the upper P-positions. Experimental results from \cite{Lar12} indicate that the upper P-positions tend to lie on two lines of slope $1.477\ldots$ and $2.247\ldots$.}\label{F1}
\end{figure}

We let $\N$ denote the positive integers and $\n$ the nonnegative integers. 
The game of 2-pile Nim is played on two piles of a finite number of tokens. Thus, its positions are represented by ordered pairs of the form  $(x,y)\in \n\times \n$. A legal move is of the form, remove a number of tokens from precisely one of the piles, at least one token and at most the whole pile. That is the set of options from any position $(x,y)$ is $$\text{Nim}(x,y)=\{(x-t,y)\mid x-t\ge 0\}\cup \{(x,y-t)\mid y-t\ge 0\}.$$ It is easy to see that the P-positions of this game are those where the pile heights are equal, that is the positions $(x, x)$, for $x\in \n$, \cite{Bou02}. We regard these positions as an infinite P-\emph{beam} of slope 1, with its source at the origin. See Figures \ref{F0} and \ref{F1}. 

In the game of Wythoff Nim a player may move as in Nim and also remove the same number of tokens from each pile, at most a whole pile, thus the set of options  from any position $(x,y)$ is $$\text{WN}(x,y)=\text{Nim}(x,y)\cup \{(x-t,y-t)\mid x-t\ge 0, y-t\ge 0\}.$$ Let $$\phi = \frac{\sqrt{5}+1}{2}$$ denote the \emph{Golden ratio}. It is known \cite{Wyt07} that a position of this game is P if and only if it belongs to the set $$\{(\lfloor \phi x\rfloor ,\lfloor\phi^2x\rfloor ), (\lfloor \phi^2 x\rfloor ,\lfloor\phi x\rfloor )\mid x\in \n\}.$$ See also Table \ref{F:1}. Thus, in the transformation from 2-pile Nim to Wythoff Nim, the single Nim-beam of P-positions has \emph{split} into two distinct beams, with sources at the origin, of slopes  $1/\phi$ and $\phi$ respectively. The intuitive meaning of the term \emph{split} is that there is an infinite sector, between two infinite regions of P-positions, which contains only N-positions.  More formally, a sequence of pairs of natural numbers $(x_i,y_i)$ \emph{splits} (or $(\alpha,\epsilon)$-\emph{splits}) if there are positive real numbers $\alpha$ and $\epsilon$ such that $$\{i>n\mid \alpha\le \frac{y_i}{x_i}\le \alpha +\epsilon\}$$ is empty for $n$ sufficiently large, but both $$\{(x_i,y_i)\mid \alpha > \frac{y_i}{x_i}\}$$ and $$\{(x_i,y_i)\mid \alpha +\epsilon < \frac{y_i}{x_i}\}$$ are infinite. 

In the game of \emph{Generalized Diagonal Wythoff Nim}, GDWN$=(p,q)$-GDWN \cite{Lar12}, in addition to the moves of Wythoff Nim, it is allowed to remove simultaneously $pt$ tokens from either of the piles and $qt$ from the other, $t\in \N$, restricted only by the number of tokens in the respective pile. That is the set of options from any position $(x,y)$ is 

\begin{align*}
(p,q)\G(x,y)=\text{WN}(x,y)&\cup\{(x-pt,y-qt)\mid x-pt\ge 0,y-qt\ge 0\}\\&\cup\{(x-qt,y-pt)\mid x-qt\ge 0,y-pt\ge 0\}.
\end{align*}

See Figure \ref{F0} for the rules of $(1,2)$-GDWN and its first few P-positions. In Figure~\ref{F1} we view the initial behavior of their respective P-beams.

Our main theorem considers the  splitting of Wythoff Nim's upper P-beam for the case $(1,2)\G$. In the last section we show that an analogous result holds for $(2,3)\G$. Before we prove these result we demonstrate a general lemma for extensions of Wythoff Nim. We show roughly that the density $\phi^{-1}$ obtained by projecting the upper P-positions of Wythoff Nim on the x-axis suffices as a lower bound for the lower asymptotic density of similar projections for any \emph{extension} of Wythoff Nim. 

\section{A natural bound for the lower asymptotic density of x-coordinates for P-positions of Wythoff Nim extensions}
An \emph{extension} of Wythoff Nim, say $G$, has a set of options $G(x,y)\supseteq$ WN$(x,y)$ from each given position $(x,y)\in \n\times\n$. The options are of the form $(x-m_1, y-m_2)\in G(x,y)$, $x\ge x-m_1\ge 0$ and $y\ge y-m_2\ge 0$, not both $m_1=m_2=0$, to disallow ties. We require that all options be symmetric, that is $(m_1,m_2)$ is a move from $(x,y)$ if and only if $(m_2,m_1)$ is a move from $(y,x)$. For a technical reason we require that at most finitely many moves of the form $(m_1,m_2)$ have the same $m_1$-coordinate (or $m_2$-coordinate), independently of from which position the move is legal. In this way the new game will have so-called \emph{complementary} sets (sequences) defining the P-positions. Two sets $\{x_i\mid i \in \N\}$ and $\{y_i\mid i\in \N\}$ are \emph{complementary} if each natural number occurs in precisely one of these sets.

It is known that the sets defining Wythoff Nim's P-positions $$A=\{\lfloor\phi n\rfloor\mid n\in \N\} \text{ and } B=\{\lfloor\phi^2 n\rfloor\mid n\in \N\}$$ are complementary, see also Figure~\ref{F:1}.  We code the \emph{upper P-positions} of a given Wythoff Nim extension $G$ by $(a_n,b_n)$ where $0<a_n< b_n$ for all $n>0$, and where the $a$-sequence is strictly increasing. By symmetry, $G$'s complete set of P-positions will be $\{(a_i,b_i),(b_i,a_i)\mid i\in \N\}\cup\{(0,0)\}$. We let $\delta_i=b_i-a_i$ for all $i$. 

\begin{Prop}\label{prop}
Any Wythoff Nim extension satisfies the following \emph{Property~W}: $a_{n-1} < a_n$, $a_n < b_n$ for all $n>0$, the sets $a = \{a_i\mid i\in \N\}$ and $b = \{b_i\mid i\in \N\}$ are complementary and, for all $i\ne j$, $\delta_i\ne \delta_j$.
\end{Prop}

\noindent{\bf Proof.} There will be no two P-positions on the same x-coordinate (or y-coordinate) since Nim-type moves are legal. It follows that it suffices that the $a$-sequence is increasing. There will be no two P-positions on any diagonal of slope 1 since the Wythoff Nim-type diagonal moves are legal. It follows that it suffices that $b_i > a_i$ for all $i > 0$. Further, for each given $x\in \N$, since there are only finitely many moves available for each $m_1$-coordinate, there has to be a least $y$ such that $(x, y)$ is a P-position. Namely, by the Nim-type condition, there is at most one P-position for each given x-coordinate less than $x$. Hence the total number of options are less than $(\sum_{i=1}^x C_i)x$, where $C_i$ denotes the finite number of moves for the $m_1$-coordinate $i$. Also, the P-positions are symmetric since the moves are symmetric. \hfill $\Box$\\

For the other direction, it is not immediately clear that any pair of sequences satisfying Property~W constitute the P-positions of some Wythoff Nim extension, but we suspect this to hold.
	
\begin{table} [ht!]
\vspace{0.5 cm}
\begin{center}
\begin{tikzpicture} [scale = 0.6]
\foreach \x/\y in {0/2, 1/2, 2/2, 3/2, 4/2, 5/2, 6/2}
 \draw[thick,brown] (\x, \y) rectangle (\x+1, \y+1);
\foreach \x/\y in {0/0, 1/0, 2/0, 3/0, 4/0, 5/0, 6/0, 0/1, 1/1, 2/1, 3/1, 4/1, 5/1, 6/1}
 \draw[thick,blue] (\x, \y) rectangle (\x+1, \y+1);
\draw (-0.5 , 1.5) node {$B_n$};
\draw (-0.5 , 0.5) node {$A_n$};
\draw (-0.5 , -0.5) node {$n$};
\draw (-0.5 , 2.5) node {$\Delta_n$};
\draw (0.5, -0.5) node {$0$};
\draw (1.5, -0.5) node {$1$};
\draw (2.5, -0.5) node {$2$};
\draw (3.5, -0.5) node {$3$};
\draw (4.5, -0.5) node {$4$};
\draw (5.5, -0.5) node {$5$};
\draw (6.5, -0.5) node {$6$};
\draw (7.5, -0.5) node {$\ldots$};
\draw (7.5, 0.5) node {$\ldots$};
\draw (7.5, 1.5) node {$\ldots$};
\draw (7.5, 2.5) node {$\ldots$};
\draw (0.5, 0.5) node {$0$};
\draw (1.5, 0.5) node {$1$};
\draw (2.5, 0.5) node {$3$};
\draw (3.5, 0.5) node {$4$};
\draw (4.5, 0.5) node {$6$};
\draw (5.5, 0.5) node {$8$};
\draw (6.5, 0.5) node {$9$};
\draw (0.5, 1.5) node {$0$};
\draw (1.5, 1.5) node {$2$};
\draw (2.5, 1.5) node {$5$};
\draw (3.5, 1.5) node {$7$};
\draw (4.5, 1.5) node {$10$};
\draw (5.5, 1.5) node {$13$};
\draw (6.5, 1.5) node {$15$};
\draw (0.5, 2.5) node {$0$};
\draw (1.5, 2.5) node {$1$};
\draw (2.5, 2.5) node {$2$};
\draw (3.5, 2.5) node {$3$};
\draw (4.5, 2.5) node {$4$};
\draw (5.5, 2.5) node {$5$};
\draw (6.5, 2.5) node {$6$};

\end{tikzpicture}
\end{center}

\caption{Wythoff Nim's upper P-positions are $(0,0),(1,2),\ldots, (A_n,B_n),\ldots$. The consecutive differences $B_n-A_n = \Delta_n$ represent the natural numbers in strictly increasing order, that is $\Delta_n = n$ for all $n$.}\label{F:1}
\end{table}

\begin{table} [ht!]
\vspace{0.1 cm}
\begin{center}
\begin{tikzpicture} [scale = 0.6]
\foreach \x/\y in {0/2, 1/2, 2/2, 3/2, 4/2, 5/2, 6/2}
 \draw[thick,brown] (\x, \y) rectangle (\x+1, \y+1);
\foreach \x/\y in {0/0, 1/0, 2/0, 3/0, 4/0, 5/0, 6/0, 0/1, 1/1, 2/1, 3/1, 4/1, 5/1, 6/1}
 \draw[thick,blue] (\x, \y) rectangle (\x+1, \y+1);
\draw (-0.5 , 1.5) node {$b_n$};
\draw (-0.5 , 0.5) node {$a_n$};
\draw (-0.5 , -0.5) node {$n$};
\draw (-0.5 , 2.5) node {$\delta_n$};
\draw (0.5, -0.5) node {$0$};
\draw (1.5, -0.5) node {$1$};
\draw (2.5, -0.5) node {$2$};
\draw (3.5, -0.5) node {$3$};
\draw (4.5, -0.5) node {$4$};
\draw (5.5, -0.5) node {$5$};
\draw (6.5, -0.5) node {$6$};
\draw (7.5, -0.5) node {$\ldots$};
\draw (7.5, 0.5) node {$\ldots$};
\draw (7.5, 1.5) node {$\ldots$};
\draw (7.5, 2.5) node {$\ldots$};
\draw (0.5, 0.5) node {$0$};
\draw (1.5, 0.5) node {$1$};
\draw (2.5, 0.5) node {$2$};
\draw (3.5, 0.5) node {$4$};
\draw (4.5, 0.5) node {$7$};
\draw (5.5, 0.5) node {$8$};
\draw (6.5, 0.5) node {$9$};
\draw (0.5, 1.5) node {$0$};
\draw (1.5, 1.5) node {$3$};
\draw (2.5, 1.5) node {$6$};
\draw (3.5, 1.5) node {$5$};
\draw (4.5, 1.5) node {$10$};
\draw (5.5, 1.5) node {$14$};
\draw (6.5, 1.5) node {$17$};
\draw (0.5, 2.5) node {$0$};
\draw (1.5, 2.5) node {$2$};
\draw (2.5, 2.5) node {$4$};
\draw (3.5, 2.5) node {$1$};
\draw (4.5, 2.5) node {$3$};
\draw (5.5, 2.5) node {$6$};
\draw (6.5, 2.5) node {$8$};

\end{tikzpicture}
\end{center}

\caption{These sequences constitute the first few P-positions of the Wythoff Nim extension $(1,2)$\G.}\label{F:2}
\end{table}

\begin{table} [ht!]
\vspace{0.4 cm}
\begin{center}
\begin{tikzpicture} [scale = 0.6]
\foreach \x/\y in {0/2, 1/2, 2/2, 3/2, 4/2, 5/2, 6/2}
 \draw[thick,brown] (\x, \y) rectangle (\x+1, \y+1);
\foreach \x/\y in {0/0, 1/0, 2/0, 3/0, 4/0, 5/0, 6/0, 0/1, 1/1, 2/1, 3/1, 4/1, 5/1, 6/1}
 \draw[thick,blue] (\x, \y) rectangle (\x+1, \y+1);
\draw (-0.5 , 1.5) node {$b_n$};
\draw (-0.5 , 0.5) node {$a_n$};
\draw (-0.5 , -0.5) node {$n$};
\draw (-0.5 , 2.5) node {$\delta_n$};
\draw (0.5, -0.5) node {$0$};
\draw (1.5, -0.5) node {$1$};
\draw (2.5, -0.5) node {$2$};
\draw (3.5, -0.5) node {$3$};
\draw (4.5, -0.5) node {$4$};
\draw (5.5, -0.5) node {$5$};
\draw (6.5, -0.5) node {$6$};
\draw (7.5, -0.5) node {$\ldots$};
\draw (7.5, 0.5) node {$\ldots$};
\draw (7.5, 1.5) node {$\ldots$};
\draw (7.5, 2.5) node {$\ldots$};
\draw (0.5, 0.5) node {$0$};
\draw (1.5, 0.5) node {$1$};
\draw (2.5, 0.5) node {$3$};
\draw (3.5, 0.5) node {$4$};
\draw[green] (4.5, 0.5) node {$7$};
\draw[blue, opacity=0.5] (4.5, 0.5) node {$7$};
\draw (5.5, 0.5) node {$8$};
\draw (6.5, 0.5) node {$9$};
\draw (0.5, 1.5) node {$0$};
\draw (1.5, 1.5) node {$2$};
\draw (2.5, 1.5) node {$5$};
\draw (3.5, 1.5) node {$6$};
\draw (4.5, 1.5) node {$10$};
\draw (5.5, 1.5) node {$13$};
\draw (6.5, 1.5) node {$15$};
\draw (0.5, 2.5) node {$0$};
\draw (1.5, 2.5) node {$1$};
\draw[red] (2.5, 2.5) node {$2$};
\draw[red] (3.5, 2.5) node {$2$};
\draw (4.5, 2.5) node {$3$};
\draw (5.5, 2.5) node {$5$};
\draw (6.5, 2.5) node {$6$};

\end{tikzpicture}
\end{center}

\caption{These sequences cannot be the initial P-positions of Wythoff Nim extension. Note that there is a partial sum of the lower sequence which is greater than the corresponding partial sum for Wythoff Nim and that this forces two differences of coordinates to coincide, since we require that the lower sequence be increasing.}\label{F:3}
\end{table}

In Tables \ref{F:2} and \ref{F:3} we show a legal and an illegal Wythoff Nim extension respectively. How do the x-coordinates of upper P-positions distribute asymptoticly? The following Lemma addresses this issue and conveys that the density $\phi^{-1}$ obtained from Wythoff Nim suffices as the \emph{lower asymptotic density} for any extension of Wythoff Nim.

\clearpage
\begin{Lemma}\label{L:2}
Suppose that the sequences $(x_i)$ and $(y_i)$ satisfy Property~W as defined in Proposition \ref{prop}. Then, for all $n\in \N$, 
\begin{align}\label{sumAsumx}
\sum_{i=0}^{n} A_i\ge \sum_{i=0}^{n} x_i
\end{align}
 and 
\begin{align}\label{sumBsumy}
\sum_{i=0}^{n} B_i\le \sum_{i=0}^{n} y_i. 
\end{align}
This implies 
\begin{align}\label{ain}
\liminf_{n\in \N}\frac{\#\{i> 0\mid x_i < n\}}{n} \ge \phi^{-1}
\end{align}
and 
\begin{align*}
\limsup_{n\in \N}\frac{\#\{i> 0\mid y_i < n\}}{n} \le \phi^{-2}.
\end{align*} 
In particular the result holds for $\{(x_i,y_i)\}$ representing the upper P-positions of any Wythoff Nim extension.
\end{Lemma}

\noindent{\bf Proof.} Denote $s(n)=\sum_{i=0}^{n} s_i$, for any sequence $(s_i)$ and $s_t(n) = \sum_{i=0}^{n} s_{t_i}$, for any sequences $(t_i)$ and $(s_{t_i})$. For a contradiction, fix an $x$-sequence and suppose that $n$ is the least number such that 
\begin{align}\label{assume1}
x(n) > A(n). 
\end{align}
Then, by the first part of property W, $A_n=x_i$ for some $i<n$ or $A_n=y_i$ for some $i<n$. Let $N=x_n-1$. Then, by minimality of $n$, $$r:=N-A_n\ge 0.$$ Define $\xi$ by $A_\xi = \lfloor \phi^{-1}A_n \rfloor$, that is $\xi = \lceil \phi^{-1}A_\xi \rceil =  \lceil \phi^{-1}\lfloor \phi^{-1}A_n \rfloor \rceil$. Then $B_{\xi} < A_n < B_{\xi+1}$ and $N = n+\xi+r$. The $N$ least numbers in $\{x_1,\ldots ,x_{n-1}\}\cup\{y_1,\ldots , y_{n-1}\}$ are 
\begin{align*}
S_1&:=\{x_1,\ldots ,x_{n-1}\}\cup\{y_{t_1},\ldots , y_{t_{\xi+r+1}}\}\\&=\{1,\dots ,N\}, 
\end{align*}
for some sequence $(t_i)$. This follows by definition of $r$, since $x$ is increasing, $x$ and $y$ are complementary and since $x_i < y_i$ (for all $i\ge n$). That is $r$ counts the number of numbers in $S_1$ greater than $A_n$. On the other hand, the $N$ least numbers in $\{A_1,\ldots ,A_{n-1}\}\cup\{B_1,\ldots , B_{n-1}\}$ are 
\begin{align*}
S_2:=\{1,\dots ,N-r-1, B_{\xi +1},\ldots , B_{\xi+r+1}\}, 
\end{align*}
where, by known properties of Wythoff's sequences, $B_{\xi+r+1}\ge N+r+1$. Observe that, by (\ref{assume1}) and definition of $r$,
\begin{align}\label{xAr}
x(n-1) \ge A(n-1)-r.
\end{align}
This implies 
\begin{align}\label{ytBr}
y_t(\xi+r+1) < B(\xi+r+1)-r.
\end{align}
 Further, by property W, we require 
 \begin{align*} 
 \delta_t (\xi+r+1)&=y_t(\xi+r+1)-x_t(\xi+r+1)\\
 &\ge \Delta(\xi+r+1)\\
 &=B(\xi+r+1)-A(\xi+r+1)\\
 &=(\xi+r+1)(\xi+r+2)/2. 
 \end{align*} 
 Hence (\ref{ytBr}) implies 
\begin{align}\label{xtAr}
x(\xi +r+1)\le x_t(\xi+r+1) < A(\xi+r+1)-r.
\end{align}
Thus, combining (\ref{xAr}) and (\ref{xtAr}), we get 
\begin{align}\label{xxAA}
x(n-1) - x(\xi+r+1) > A(n-1) - A(\xi+r+1).
\end{align}
We make the following partitioning of the set $S_1$, 
\begin{align*}
X^- &= \{1,\ldots , x_{\xi+r+1}\}\\
X^+ &= \{ x_{\xi+r+2},\ldots , x_{n-1}\}\\
Y^- &= \{y_i\mid y_i \le A_{\xi+r+1}\},\\
Y^+ &= \{y_i \mid A_{\xi+r+1} < y_i\le N\}
\end{align*}
and put $Y = Y^-\cup Y^+$. See also Figures \ref{F:4} and \ref{F:5}. In other words $\sum_{s\in Y} s = y_t(\xi+r+1)$ and, as we have seen in (\ref{ytBr}), this sum is \emph{small} as compared to the corresponding sum on the $B$-sequence. But the inequality in (\ref{ytBr}) is obtained via (\ref{xAr}) and by comparing the sets $S_1$ and $S_2$. In particular, there is a bias towards \emph{smaller} numbers in $X^-$ ($\le A_{\xi+r+1}$) in comparison with the numbers in the $A$-sequence $\le A_{\xi+r+1}$. 
Therefore, again by (\ref{assume1}), there are relatively many large numbers in $X^+$ in comparison with the $A$-sequence. Again, by complementarity of our sequences, we conclude that there is a bias towards \emph{smaller} numbers in $Y^+$ than what was conveyed by (\ref{ytBr}). Hence, since we are only considering integer entries, (\ref{xxAA}) implies that we can strengthen  (\ref{ytBr}) to 
\begin{align}\label{ytBr1}
y_t(\xi+r+1) < B(\xi+r+1)-r-1.
\end{align}
But then we can repeat all arguments until (\ref{xxAA}) which this time rather becomes 
\begin{align*}
x(n-1) - x(\xi+r+1) > A(n-1) - A(\xi+r+1) + 1,
\end{align*}
which implies 
\begin{align*}
y_t(\xi+r+1) < B(\xi+r+1)-r-2,
\end{align*}
and so on. Since the left hand side cannot be arbitrarily small, but rather $y_t(\xi+r+1)\ge (\xi+r+1)(\xi+r+2)/2$, our argument gives a contradiction.
\begin{figure} [ht!]
\vspace{0.5 cm}
\begin{center}
\begin{tikzpicture} [scale = 1]
\foreach \x/\y in {0/1}
 \draw[thick,blue] (\x, \y) rectangle (\x+5, \y+1);
\foreach \x/\y in {0/0}
 \draw[thick,blue] (\x, \y) rectangle (\x+5, \y+1);
\foreach \x/\y in {5/0}
 \draw[thick,blue] (\x, \y) rectangle (\x+3, \y+1);
\draw (0.6 , .5) node {$X^-$};
\draw (0.5 , 1.5) node {$Y$};
\draw (5.6 , .5) node {$X^+$};
\draw (4.1, 0.4) node {$x_{\xi+r+1}$};
\draw (4.1, 1.4) node {$y_{\xi+r+1}$};
\draw (7.3, 0.4) node {$x_{n-1}$};
\draw[red] (6.8, 2.9) node {$A_n$};
\draw[thick, red, dashed, ->] (6, 3) ..controls (4.6, 3.1) and (3.5, 2.5) .. (3, 1.68);
\draw[thick, red, dashed, ->] (7.2, 2.3) ..controls (7.5, 1.5) and (7.2, 1) .. (7.1, .76);
\end{tikzpicture}
\end{center}

\caption{The region $X^-$ contains the $\xi+r+1$ least numbers in the $x$-sequence, whereas $Y$ contains the $\xi+r+1$ least numbers in $y$. In this figure they coincide with the first entries, but this is usually not the case, since the $y$-sequence is not necessarily increasing. The region $X^+$ contains the entries $x_{\xi+r+2},\ldots , x_{n-1}$. We have that $X^-\cup X^+\cup Y = \{1,\ldots ,N\}$ and  $(X^-\cup X^+)\cap Y = \emptyset$. Wythoff Nim's entry $A_n$ belongs either to $Y$ or $X^+$. This fact, which is exploited in the proof, makes it too crowded in some region.}\label{F:4}
\end{figure}
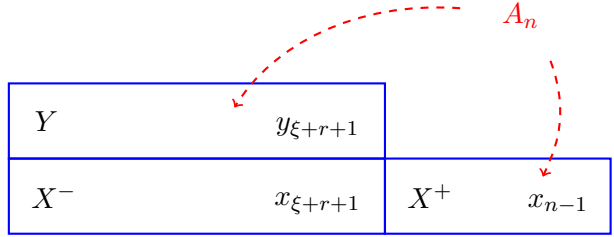
	
\begin{figure}[ht!]
\vspace{0.5 cm}
\begin{center}
\begin{tikzpicture} [scale = 1]
\foreach \x/\y in {0/1}
 \draw[thick,blue] (\x, \y) rectangle (\x+5.2, \y+1);
\foreach \x/\y in {5.2/1}
 \draw[thick,blue] (\x, \y) rectangle (\x+3.8, \y+1);
\foreach \x/\y in {0/0}
 \draw[thick,blue] (\x, \y) rectangle (\x+5, \y+1);
\foreach \x/\y in {5/0}
 \draw[thick,blue] (\x, \y) rectangle (\x+4, \y+1);
\draw (0.6 , .5) node {$X^-$};
\draw (0.5 , 1.5) node {$Y^-$};
\draw (5.6 , .5) node {$X^+$};
\draw (5.8 , 1.5) node {$Y^+$};
\draw (4.1, 0.4) node {$x_{\xi+r+1}$};
\draw (4.2, 1.5) node {$A_{\xi+r+1}$};
\draw (8.3, 0.4) node {$x_{n-1}$};
\draw (8.5, 1.5) node {$N$};
\draw[red] (7.1, 1.4) node {\small{smaller}};
\draw[red] (2, .4) node {\small{small}};
\draw[blue] (2, 1.4) node {\small{large}};
\draw[blue] (7, .4) node {\small{large}};
\end{tikzpicture}
\end{center}

\caption{The region $Y^-$ contains the numbers in $Y$ less than $ A_{\xi +r + 1}$ and $Y^+=Y\setminus Y^-$. }\label{F:5}
\end{figure}

Hence (\ref{sumAsumx}) holds and by complementarity (\ref{sumBsumy}). Note also that, by (\ref{sumAsumx}), for all $n$,
\begin{align}\label{nnxph}
n(n+1)/2\ge x(n)\phi^{-1},
\end{align}
which implies $x_n\le \phi n+o(n)$. Hence $$\frac{\#\{i\mid x_i\le x_n \}}{x_n}=\frac{n}{x_n} \ge \frac{1}{\phi +o(1)},$$ 
which implies
\begin{align*}
\liminf_{n\in \N}\frac{\#\{i> 0\mid x_i < n\}}{n} \ge \phi^{-1}.
\end{align*}
\hfill $\Box$

\section{The splitting of $(1,2)\G$}
In this section we analyze the game $(1,2)\G$ from \cite{Lar12} and prove that its upper P-positions $(2,0.05)$-split. 
The following Lemma shows that it suffices to establish a positive lower asymptotic density of x-coordinates of P-positions above the line $y=2x$.
\begin{Lemma}\label{L:3}
If there is a positive lower asymptotic density of x-coordinates of P-positions above the line $y=2x$, then the upper P-positions $\{(a_n,b_n)\mid n\in \n\}$ of $(1,2)\G$ split.
\end{Lemma}
\noindent{\bf Proof.}
Let $(k_i)$ denote the unique increasing sequence of indices such that, for all $0<i $, $b_{k_i}/a_{k_i} > 2$ and if $k_i<j<k_{i+1}$ then $b_j/a_i<2$ (put $k_0=0$).
Already in \cite{Lar12} we proved that there are infinitely many such $i$ and that the $b$-sequence satisfies, for all $i$, 
\begin{align}\label{+1}
b_{k_{i+1}} = 2(a_{k_{i+1}} - a_{k_{i}})+b_{k_{i}}+1. 
\end{align}
The assumption is that there is an $\epsilon >0$ such that, for all sufficiently large $N$, 
\begin{align}\label{iN}
\#\{i\mid a_{k_i}<N\}\ge \epsilon N. 
\end{align}
It suffices to show that, for all sufficiently large $N$, the set $$\{(N, 2N),(N, 2N+1),\ldots ,(N, (2+\epsilon)N\}$$ will contain only N-positions (assuming that $\epsilon N$ is an integer). Suppose that there is a move from $(N,2N+x)$ to $(a_{k_i},b_{k_i})$. Then there is a $t$ such that $2N+x-b_{k_i}=2t$ and $N-a_{k_i}=t$. By (\ref{iN}) it suffices to show that there is a move from $(N,2N+x+1)$ to $(a_{k_{i+1}},b_{k_{i+1}})$, because there is a move from $(N,2N)$ to $(0,0)$. By eliminating $t$, we get $2N+x+1-b_{k_i}=2(N-a_{k_i})+1$ which, by (\ref{+1}) becomes $2N+x+1-b_{k_{i+1}}=2(N-a_{k_{i+1}})$ and we are done with this part.

 By \cite{Lar12} we know that there are infinitely many P-positions $(a_n,b_n)$ satisfying $1\le\frac{b_n}{a_n}\le 2$. Hence, given the assumptions, each requirement for a split is satisfied.
 \hfill $\Box$\\

\begin{Thm}
The upper P-positions $\{(a_n,b_n)\mid n\in \n\}$ of $(1,2)\G$ split. 
\end{Thm}
\noindent{\bf Proof.} If there is a positive lower density of P-positions above the line $y=2x$ then, by Lemma \ref{L:3} the result holds.

Assume for a contradiction that there are at most $o(N)$ P-positions above the line $y=2x$ with x-coordinate less than $N$. Then almost all upper P-positions lie in the sector defined by the lines $y=x$ and $y=2x$. By Lemma~\ref{L:2} the number of upper P-positions with x-coordinate less than $N$ is at least $\phi^{-1}N - o(N)$.

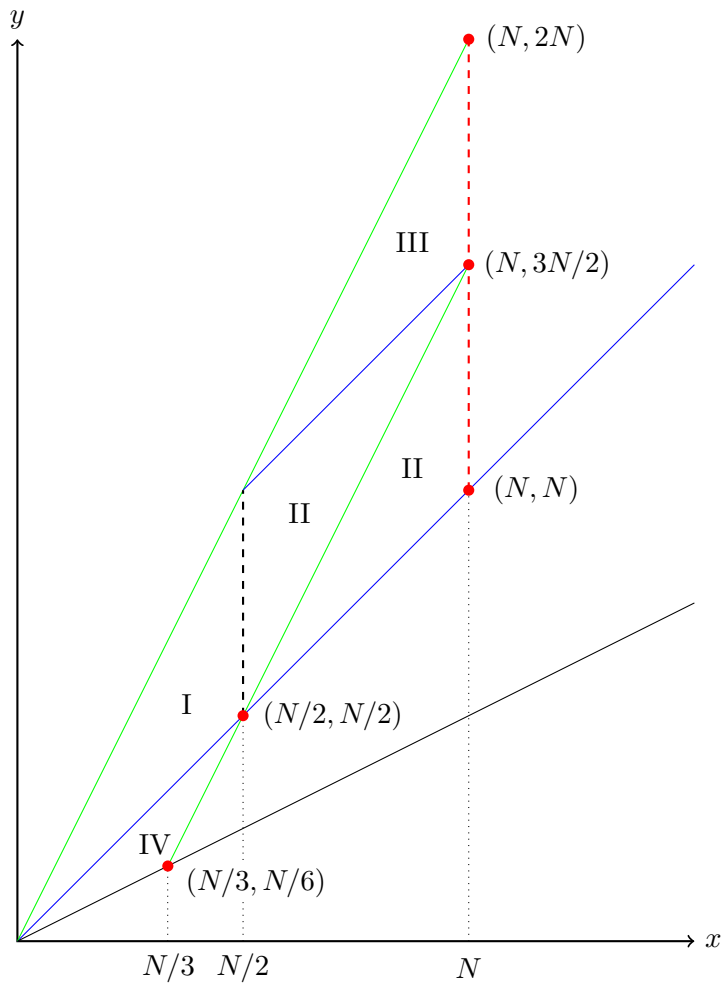
\begin{figure}[ht!]
\begin{center}
\begin{tikzpicture}[scale = 3]
    \draw [<->,thick] (0,4) node (yaxis) [above] {$y$}
        |- (3,0) node (xaxis) [right] {$x$};
    \draw[blue] (0,0) coordinate (a_1) -- (3,3) coordinate (a_2);
 \draw[green] (0,0) -- (2,4);
 \draw[green] (2/3,1/3) -- (2,3);
 \draw (0,0) -- (3,3/2);
 \draw[blue] (1,2) -- (2,3);
 \draw[dashed,red,thick] (2,2) -- (2,4);
 \draw[dotted] (1,0) -- (1,.18);
 \draw[dotted] (1,0.36) -- (1,1);
\draw[dotted] (.666,0) -- (.666,.333);
 \draw[dashed,thick] (1,1) -- (1,2);
 \draw[dotted] (2,0) -- (2,2);
\draw (0.75, 1.05) node {I};
\draw (1.75, 2.1) node {II};
\draw (1.25, 1.9) node {II};
\draw (1.75, 3.1) node {III};
\draw (0.61, .43) node {IV};
\draw (2,-0.12) node {$N$};
\draw (1,-0.12) node {$N/2$};
\draw (.67,-0.12) node {$N/3$};
\draw (2.3,2) node {$(N,N)$};
\draw (2.35,3) node {$(N,3N/2)$};
\draw (2.3,4) node {$(N,2N)$};
\draw (1.4,1) node {$(N/2,N/2)$};
\draw (1.06,.26) node {$(N/3,N/6)$};
\fill[red] (2/3,1/3) circle (.7pt);    
\fill[red] (2,2) circle (.7pt);    
\fill[red] (1,1) circle (.7pt);    
\fill[red] (2,3) circle (.7pt);    
\fill[red] (2,4) circle (.7pt);
\end{tikzpicture}\caption{We consider P-positions in four regions sending beams of N-positions to lattice points with x-coordinate $N$. The ones from region I and II have slope 1 and the ones from region I, III and IV have slope 2. The y-coordinate 3N/2 is pivotal in our argument. Below it we count the beams of slope 1 and above it the beams of slope 2.}\label{fig5}
\end{center}
\end{figure}

We get $$\tau = \tau(N):=\frac{\#\{a_i\mid i\le N\}}{N}\ge \phi^{-1},$$ for all sufficiently large $N$. %, where the last inequality is by Lemma \ref{L:2}. 
By this lemma, the number of N-beams of slope 1, which intersect the red dotted line in Figure \ref{fig5} between $y=N$ and $y=\frac{3N}{2}$, originating from P-positions in region
\begin{itemize} 
\item [(I)] is at least $\tau\frac{N}{2}$,
\item [(II)] is at least $c\tau\frac{N}{2}$, where $0<c<1$ is some real constant.
\end{itemize}
Similarly, the number of N-beams of slope 2, which intersects the red dotted line between $y=\frac{3N}{2}$ and $y=2N$, originating from P-positions in region
\begin{itemize} 
\item [(I)] is at least $\tau\frac{N}{2}$,
\item [(III)] is at least $(1-c)\tau\frac{N}{2}$, where $0<c<1$ is the real constant from (II),
\item[(IV)] is at least $(1-\tau)\frac{N}{3}$,
\end{itemize}
where the last item is by complementarity of the $a$ and $b$ sequences.
Hence, the number of N-beams of slope 1, from regions I and II, is at least $$(1+c)\frac{\tau N}{2}.$$ The number of N-beams of slope 2, from regions I, III and IV, is at least $$\left(\frac{\tau}{2}+(1-c)\frac{\tau}{2} + (1-\tau)\frac{1}{3}\right)N.$$ From these expressions we see that the constant $c$ does not affect the total number of N-beams. Hence we omit it and, by $\tau\ge \phi^{-1}$, get that the total number of N-beams is at least $N\frac{7\phi^{-1}+2}{6}>1.05N$, which contradicts the definition of N and P, namely it implies that either there are two P-positions on the same line of slope 2 or there are two P-positions on the same line of slope 1. \hfill $\Box$\\

\clearpage

\section{The game $(2,3)\G$}
By an analogous method one can also prove that the upper P-positions of $(2,3)\G$ split. Some extra care is needed for the treatment of the P-positions above the line of slope 3/2, but one can see that the packing of the N-beams originating from P-positions above this line will be dense also for this game, within an $O(N)$ distance to the y-coordinate $3N/2$ (with $N$ even), although they will not be strictly ``greedy" as for $(1,2)\G$. Another slight complication is that one needs to regard two columns, $N, N+1$ for this game, rather than the single $N$-column in $(1,2)\G$. However this still implies that the number of positions to check for its N-status between this line and the one of slope 1, is $N$. Further, the contribution from analogues of regions I,II and III in Figure \ref{fig5} give the same estimate as for $(1,2)\G$. By Lemma \ref{L:2} and, by inspection, for this case the contribution from region IV gives an additional number of P-positions below x-coordinate $3N/10$, at least $3(1-\tau)N/10$ of them. Hence the total number is at least $(\frac{6\tau}{5} + \frac{3}{10})N>1.04N$.

\begin{Thm}
The upper P-positions $\{(a_n,b_n)\mid n\in \n\}$ of $(2,3)\G$ split. 
\end{Thm}

For other variations of GDWN the analysis seems more technical and new ideas may be needed.
\section{Questions} Property W as defined in Proposition \ref{prop} is in fact a property of a given set of positive integers as follows. Let $\{s_i\}=S\subset \N$ (with the $s_i$'s distinct). Then we have demonstrated that, if there is an ordering of the numbers in $\{t_i\} = \N\setminus S$ such that $t_i - s_i = t_j - s_j$ implies $i=j$, then the lower asymptotic density of $S$ must be greater than or equal to $\phi^{-1}$. 

Does the converse hold? That is, if there is no such ordering of the entries in the complement of $S$, is the lower density necessarily strictly less than  $\phi^{-1}$? For a counterexample find a set $S$ without this property which has lower density greater than or equal to $\phi^{-1}$. 

A sequence $(x_i,y_i)$ \emph{density-splits} if it $(\alpha,\epsilon)$-splits and each of the sets $\{x_i\mid \alpha\ge \frac{y_i}{x_i}\}$ and $\{x_i\mid \alpha +\epsilon \le \frac{y_i}{x_i}\}$ has a positive lower asymptotic density. Do our results hold if we exchange split for density-split everywhere? We conjecture a positive answer, but some details are still missing for the upper P-positions' lower P-beam, see \cite{Lar12}.

\end{document}